\documentclass[a4paper]{article}

\usepackage[latin1]{inputenc}
\usepackage[french,british]{babel}
\usepackage[T1]{fontenc}
\usepackage{setspace}
\usepackage{indentfirst} 
\usepackage{fancyhdr} 
\usepackage{nicefrac}
\usepackage{textcomp}
\usepackage{lastpage}
\usepackage{amssymb,amsfonts,amsmath,amsthm}
\usepackage{latexsym}
\usepackage{url}
\usepackage{verbatim}
\usepackage[all]{xy}
\usepackage{mathrsfs}
\usepackage{stmaryrd}
\usepackage{oldgerm}
\usepackage{bm}
\usepackage{graphics}
\usepackage{wrapfig}
\usepackage{fancybox}
\usepackage{graphicx}
\input xy
\xyoption{all}

\theoremstyle{plain}
\newtheorem{theo}{Theorem}[section]
\newtheorem*{theo*}{Theorem}
\newtheorem{lem}[theo]{Lemma}
\newtheorem*{lem*}{Lemma}

\newtheorem{prop}[theo]{Proposition}
\newtheorem*{prop*}{Proposition}

\newtheorem*{conj*}{Conjecture}

\theoremstyle{remark}
\newtheorem{rem}[theo]{Remark}

\theoremstyle{definition}
\newtheorem{defi}[theo]{Definition}

\numberwithin{equation}{section}

\DeclareMathAlphabet{\mathpzc}{OT1}{pzc}{m}{it}
\DeclareMathOperator{\NN}{\mathbb{N}}

\DeclareMathOperator{\QQ}{\mathbb{Q}}

\DeclareMathOperator{\HH}{\mathbb{H}}

\DeclareMathOperator{\Hom}{Hom}

\DeclareMathOperator{\OO}{\mathscr{O}}

\DeclareMathOperator{\M}{\mathbf{M}}

\DeclareMathOperator{\Spec}{Spec}

\DeclareMathOperator{\h}{H}

\DeclareMathOperator{\rank}{rank}

\DeclareMathOperator{\Frac}{Frac}

\DeclareMathOperator{\f}{F}

\DeclareMathOperator{\Isoc}{Isoc}

\DeclareMathOperator{\Con}{Con}

\begin{document}

\selectlanguage{british}

\thispagestyle{empty}

\title{Full faithfulness for overconvergent $\f$-de Rham--Witt connections\\
Fid\'elit\'e pleine pour les $\f$-connexions de de Rham--Witt surconvergentes}
\author{Ertl Veronika}
\date{}
\maketitle

\begin{abstract}
{\noindent
Let $X$ be a smooth variety over a perfect field of characteristic $p>0$. We define overconvergent $F$-de Rham-Witt connections as an analogue for $F$-crystals over proper schemes. We prove that the forgetful functor from the category of overconvergent $F$-de Rham--Witt connections to the category of convergent $F$-de Rham--Witt connections is fully faithful.}
\end{abstract}

\vspace{-9pt}

\selectlanguage{french}

\begin{abstract}
{\noindent
Soit $X$ une vari\'et\'e lisse sur un corps parfait de charact\'eristique  $p>0$. On d\'efinit les $F$-connections de de Rham--Witt surconvergentes comme analogues de $F$-cristaux sur les sch\'emas propre. On montre que le foncteur oubli de la cat\'egorie de $F$-connections de de Rham--Witt surconvergentes \`a la cat\'egorie de $F$-connections de de Rham--Witt convergentes est pleinement fid\`ele. 
}
\end{abstract}

\selectlanguage{british}


\let\thefootnote\relax\footnote{\textit{Scientific field:} Number theory}

\addcontentsline{toc}{section}{Introduction}
\section*{Introduction}

Let $X$ be a smooth variety over a perfect field $k$ of positive characteristic $p$. In \cite{Bloch3} Spencer Bloch studies connections on the de Rham--Witt complex of $X$ and establishes the following equivalence of categories.

\begin{theo*}[Bloch 2001, \protect{\cite[Cor. 1.2]{Bloch3}}]  The functor $\M\mapsto(M,\nabla)$ defines an equivalence of categories between locally free $F$-crystals on $X$ and locally free $W(\OO_X)$-modules with quasi-nilpotent integrable connections and a Frobenius structure. \end{theo*}

Unfortunately, the category of (convergent) crystals and the associated connections via this functor does not seem to be the right category to consider in the case of non-proper varieties. As explained by Davis, Langer and Zink in \cite{DavisLangerZink}, there is a remedy by restricting the de Rham--Witt complex to so-called overconvergent elements, which form again a complex over the Witt vectors $W(k)$ of $k$. In \cite{Ertl2} we introduce a category of coefficients for this overconvergent de Rham--Witt complex, which we call overconvergent de Rham--Witt connections, and show that this is a reasonable category to be considered in the sense, that the resulting cohomology theory compares well to Monsky--Washnitzer cohomology with (overconvergent) coefficients and rigid cohomology with overconvergent isocrystals. 

It does not seem to be clear at all that the category of overconvergent isocrystals embeds fully faithful in the category of convergent isocrystals. Yet Kiran Kedlaya is able to prove this statement if he considers Frobenius action as an additional structure. His statement in \cite{Kedlaya4} uses Berthelot's notion of convergent and overconvergent $F$-isocrystals, denoted for some integer $a\in\NN$ by
$F^a$-$\Isoc(X/K_0)$ and $F^a$-$\Isoc^\dagger(X/K_0)$, where $K_0$ is the fraction field of $W(k)$.

\begin{theo*}[Kedlaya 2003] 
The forgetful functor
$$j^\ast:F^a\text{-}\Isoc^\dagger(X/K_0)\rightarrow F^a\text{-}\Isoc(X/K_0)$$
is fully faithful.
\end{theo*}

In this small note we attempt to obtain a similar result for de Rham--Witt connections.

In Section 1 we introduce the notion of convergent and overconvergent $F$-de Rham--Witt connections. More precisely,
\begin{enumerate}
\item the category of convergent, locally free, quasinilpotent, integrable $(\sigma,\nabla)$-modules over the Witt sheaf $W(\OO_X)$, denoted by $F$-$\Con(X/W(k))$ and 
\item the category of overconvergent, locally free, quasinilpotent, integrable $(\sigma,\nabla)$-modules over $W^\dagger(\OO_X)$, denoted by $F$-$\Con^\dagger(X/W(k))$.
\end{enumerate}

In the same section we also recall the notion of Newton slopes, as this is vital for the prove of the main theorem, namely:

\begin{theo*}[Main result]
The forgetful functor
$$j^\ast: F\text{-}\Con^\dagger(X/W(k))\rightarrow F\text{-}\Con(X/W(k))$$
is fully faithful.
\end{theo*}

In the second section we prove a local version of this, and show in Section 3 that this is in fact sufficient by a standard reduction from the global to the local setting.

\subsubsection*{Acknowledgements}

Special thanks are due to Kiran Kedlaya, who brought this question to my attention, and provided important input through his paper ``Full faithfulness for overconvergent $F$-isocrystals'' \cite{Kedlaya4}. The author is thankful to the Universit\"at of Regensburg for continued support for her research.

\section{Definitions and notations}

This paper follows up a question that came up in \cite{Ertl2}  and so we will recall and build upon its notation and terminology. Throughout this paper, let $k$ be a perfect field of characteristic $p>0$, $W(k)$ its ring of $p$-typical Witt vectors and $K_0=\Frac(W(k))$ the fraction field. Denote by $\sigma$ the absolute Frobenius automorphism on $W(k)$. Let $X$ be a smooth variety over $k$.

 For a finite $k$-algebra $\overline{A}$ let $W^\dagger(\overline{A})\subset W(\overline{A})$ be the ring of overconvergent Witt vectors as defined by Davis, Langer and Zink in \cite{DavisLangerZink2} and later extended to the overconvergent de Rham--Witt complex $W^\dagger\Omega_{\overline{A}}\subset W\Omega_{\overline{A}}$. Both can in fact be globalised to a sheaf on $X$ over $k$, denoted by $W^\dagger\Omega_X$ and  $W^\dagger(\OO_X)$ respectively. 

In order to extend the notion of a crystal on a proper variety, we defined in \cite[Definition 1.3.1]{Ertl2} an overconvergent de Rham--Witt connection on $X$ to be a $W^\dagger(\OO_X)$-module $M$ together with a map
$$\nabla: M\rightarrow M\otimes_{W^\dagger(\OO_X)} W^\dagger\Omega^1_{X/k}$$
satisfying the Leibniz rule. It is said to be integrable if for the induced map in the complex  $\nabla^2=0$ holds. This is the overconvergent analogue of the definition for a de Rham--Witt connection found in \cite{Bloch3}. A morphism of (overconvergent) de Rham--Witt connections $(M,\nabla_M)\rightarrow (N,\nabla_N)$ is a morphism of the associated complexes. We denote the category of locally free, quasinilpotent (in the sense of \cite{BerthelotOgus}), integrable de Rham--Witt connections by $\Con(X/W(k))$ and the subcategory of locally free, quasinilpotent, integrable, overconvergent de Rham--Witt connections by $\Con^\dagger(X/W(k))$.


The definition of the de Rham--Witt complex provides by functoriality a Frobenius and Verschiebung map on $W\Omega_{X/k}$ in degree $i$ divisible by $p^i$, which descend to the overconvergent subcomplex $W^\dagger\Omega_{X/k}$ (cf. \cite[Section 2.4]{Ertl}). 

\begin{defi}
A $\sigma$-module over $W(\OO_X)$ is a $W(\OO_X)$-module $M$ of finite type equipped with a $\sigma$-linear map $F:M\rightarrow M$ which becomes an isomorphism over $K_0$. A morphism $\psi:M\rightarrow N$ between two $\sigma$-modules over $W(\OO_X)$ is a morphism of $W(\OO_X)$-modules such that $F\circ \psi=\psi\circ F$. \\
We define a $(\sigma,\nabla)$-module over $W(\OO_X)$ to be a $\sigma$-module $M$ together with a connection $\nabla:M\rightarrow M\otimes W\Omega^1_{X/k}$  that makes the diagram
$$\xymatrix{M\ar[r]^{\nabla\qquad}\ar[d]_{F} & M\otimes W\Omega^1_{X/k}\ar[d]^{F\otimes F}\\
M\ar[r]^{\nabla\qquad}  & M\otimes W\Omega^1_{X/k}}$$
commute.
\end{defi}

In other words, a $(\sigma,\nabla)$-module over $W(\OO_X)$ is a de Rham--Witt connection in the sense of \cite{Bloch3}  with a $\sigma$-linear map $F$ that is horizontal for the connection $\nabla$. It is said to be integrable (resp. quasinilpotent) if $\nabla$ is integrable (resp. quasinilpotent).

\begin{defi}
Given a $(\sigma,\nabla)$-module $M$ over $W(\OO_X)$, we define it's Tate twist $M(\ell)$ as the same underlying $W(\OO_X)$-module with the action of $F$ multiplied by $p^\ell$.
\end{defi}

Bloch proves the following \cite[Corollary 1.2]{Bloch3}.

\begin{prop}
There is an equivalence of categories between the category of locally free $F$-crystals on $X$ and the category of locally free, quasinilpotent, integrable $(\sigma,\nabla)$-modules over $W(\OO_X)$.
\end{prop}

We denote the category of locally free, quasinilpotent, integrable $(\sigma,\nabla)$-modules over $W(\OO_X)$ by $F$-$\Con(X/W(k))$. 


We recall the definition of Newton slope for a crystal as introduced for example in \cite{Katz}, and which carries over to $\sigma$-modules and by extension to $(\sigma,\nabla)$-modules. In the case of a point we have the following. Let $(M,F)$ be a $\sigma$-module over $W(k)$ of rank $r$. 




There are several ways to define Newton slopes. For $\lambda=\frac{t}{s}\in\QQ$, where $t$ and $s$ are relatively prime, we denote 
$$E(\lambda)= \left(W(k)[T]/(T^s-p^t)\right)$$
the $\sigma$-module with obvious Frobenius (which is in fact injective). The rational number $\lambda$ is called the slope of $E(\lambda)$. If $k$ is algebraically closed a well known theorem of Dieudonn\'e--Manin tells us that the category of $\sigma$-modules up to isogeny is semi-simple and the simple objects are given by modules of the form $E(\lambda)$. Thus each $\sigma$-module $(M,F)$ is isogenous to a direct sum of elementary ones
$$(M,F)\sim\bigoplus_{\lambda\in\QQ_+}E(\lambda)^{n_\lambda}.$$

\begin{defi}
The $\lambda_1<\cdots<\lambda_q$ such that $n_{\lambda_i}\neq 0$ are the Newton slopes of $M$, and $n_\lambda$ is the multiplicity of $\lambda$. If $k$ is not algebraically closed, the slopes of $M$ are by definition the slopes of $M\otimes W(\overline{k})$.
\end{defi}

For the next characterisation of Newton slope we follow \cite{Katz} and choose a natural number $N\in\NN$ divisible by $r!$ where $r$ is the rank of the $\sigma$-module $(M,F)$, and consider the discrete valuation ring
$$R=W(\overline{k})[X]/(X^N-p)=W(\overline{k})[p^{\frac{1}{N}}].$$
The obvious extension of $\sigma$ to $R$ is to set $\sigma(X)=X$. Let $K$ be the fraction field of $R$. By Dieudonn\'e--Manin $M\otimes_{W(k)}K$ admits a $K$-basis $\{e_1,\ldots,e_r\}$ of Frobenius eigenvectors, i.e.
$$F(e_i)=p^{\lambda_i}e_i$$
where $\lambda_1\leqslant\cdots\leqslant\lambda_r$ are again the slopes (this time appearing with multiplicity). This says essentially that over $K$ the Frobenius morphism $F$ of $M$ is diagonalisable. 

Over $R$ we are not as lucky, but can at least say, that the Frobenius is trigonalisable. More precisely, there is an $R$-basis $\{u_1,\ldots,u_r\}$ of $M\otimes_{W(k)}R$ such that $F$ is given by an upper triangular matrix with $p^{\lambda_i}$ on the diagonal
$$\left(
\begin{array}{cccc}
p^{\lambda_1} & \alpha_{12} &\ldots& \alpha_{1r}\\
0 & p^{\lambda_2} & \ldots & \alpha_{2r}\\
\vdots &\vdots& \ddots &\vdots\\
0 & 0 &\ldots & p^{\lambda^r}\\
\end{array}
\right)$$
where $\alpha_{ij}\in R$, or in other words
$$F(u_i)\cong p^{\lambda_i}\mod\sum_{j<i}R u_j.$$

\begin{defi}
The Newton polygon of $(M,F)$ is the graph of the Newton function on $[0,r]$ given by
$$N_{(M,F)}(i)=\begin{cases} 
	0 &\quad\text{ if } i=0,\\
	\lambda_1+\cdots+\lambda_i &\quad\text{ if } i\in\{1,\ldots,r\},\\
	\text{extended linearly} & \text{ if } i\neq\NN_0.
\end{cases}$$
\end{defi}


This definition that will also serve as the core in the case when $X$ is not a point, as we talk about a local property. Although Katz gives in \cite{Katz} a much more general definition, we will stick to the case, when $X$ is a smooth scheme over the perfect field $k$, i.e. locally it is of the first type mentioned in \cite[Sec. 2.4]{Katz}. 

\begin{defi}
Let $(M,F)$ be a $\sigma$-module on $W(\OO_X)$, and $\mathfrak{p}\in X$ a point. The Newton slopes of $M$ at the point $\mathfrak{p}$ are the Newton slopes of $(M_{\kappa(\mathfrak{p})},F_{\kappa(\mathfrak{p})})$, where we localise at $\mathfrak{p}$ and $\kappa(\mathfrak{p})$ denotes the local field of $\mathfrak{p}$. 
\end{defi}


We restrict these concepts now to the category of overconvergent modules.

\begin{defi}
An overconvergent $\sigma$-module is a $W^\dagger(\OO_X)$-module $M$ of finite type with a $\sigma$-linear injective map $F:M\rightarrow M$. A morphism $\psi:M\rightarrow N$ between two $\sigma$-modules over $W^\dagger(\OO_X)$ is a morphism of $W^\dagger(\OO_X)$-modules such that $F\circ \psi=\psi\circ F$. \\
We define an overconvergent $(\sigma,\nabla)$-module over $X$ to be an overconvergent $\sigma$-module $M$ together with a connection $\nabla:M\rightarrow M\otimes W^\dagger\Omega^1_{X/k}$  such that $F$ is horizontal for $\nabla$ as in the convergent definition above.
\end{defi}

In other words, an overconvergent $(\sigma,\nabla)$-module is an overconvergent de Rham--Witt connection in the sense of \cite{Ertl2} that is compatible with a $\sigma$-linear map $F$. Integrability and quasinilpotence are the same as above. 

We denote the category of locally free, quasinilpotent, integrable, overconvergent $(\sigma,\nabla)$-modules by $F$-$\Con^\dagger(X/W(k))$.

\section{Local full faithfulness}

The goal of this section is to prove the main result locally. That is, we want to prove that locally the forgetful functor from overconvergent to convergent $F$-de Rham--Witt connections is fully faithful -- in other words, that the induced map on the associated homomorphism groups is bijective. The argumentation here is inspired by Kedlaya's in the case of (overconvergent) $F$-isocrystals in \cite{Kedlaya4}.

The setting to be considered is the following. Let $\overline{A}$ be a finite $k$-algebra and $W^\dagger(\overline{A})\subset W(\overline{A})$ and $W^\dagger\Omega_{\overline{A}/k}\subset W\Omega_{\overline{A}/k}$ as above. Let $M_1$ and $M_2$ be locally free, quasinilpotent, integrable, overconvergent $(\sigma,\nabla)$-modules over $\overline{A}$. As we are looking at the situation locally, we can even assume that $M_1$ and $M_2$ are free modules over $W^\dagger(\overline{A})$. Let 
$$j^\ast: F\text{-}\Con^\dagger(X/W(k))\rightarrow F\text{-}\Con(X/W(k))$$
be the forgetful functor. In the above situation it induces a map
$$\Hom_{F\text{-}\Con^\dagger(\overline{A}/W(k))}(M_1,M_2)\rightarrow \Hom_{F\text{-}\Con(\overline{A}/W(k))}(j^\ast M_1,j^\ast M_2).$$
It is clear that this map is injective. Therefore, our task is to prove surjectivity.

We start by proving two lemmas in order to set up the proof of the main statement of this section. 

\begin{lem}
Let $M$ be a non-zero $\sigma$-module over $W^\dagger(\overline{A})$, and let $\varphi:M\rightarrow W(\overline{A})$ be a $W^\dagger(\overline{A})$-linear injective map such that for some integer $\ell$
$$\varphi(F\mathbf{e})=p^\ell\varphi(\mathbf{e})^\sigma,$$
for all $\mathbf{e}\in M$. Then $M$ has rank one and is of slope $\ell$ at every point $\mathfrak{p}$ on $\Spec \overline{A}$. 
\end{lem}

\begin{proof} 
Let $\kappa=\kappa(\mathfrak{p})$ be the local field at the point $\mathfrak{p}$ and consider the localisation $M_{\mathfrak{p}}$ of $M$ at $\mathfrak{p}$. Let $R= W(\overline{\kappa})[X]/(X^N-p)$ where $N\in\NN$ is big enough as in Section 1, and let $K$ be the fraction field of $R$. By linearity, we extend the localisation of $\varphi$ at $\mathfrak{p}$ to the $\sigma$-module $M_{\mathfrak{p}}\otimes_{W(\overline{\kappa})}K$, which is possible since all elements of $W(\overline{\kappa})$ are overconvergent and $R$ is a finite extension of $W(\overline{\kappa})$. But $M\otimes_{W(\overline{\kappa})}K$ has a $K$-basis $\{\mathbf{e}_1,\ldots,\mathbf{e}_r\}$ that diagonalises  $F$. For each $i$, let $p^{\lambda_i}$ be the eigenvalue associated to $\mathbf{e}_i$. It follows then that
\begin{equation}\label{EquationNurnberg}p^{\lambda_i}\varphi(\mathbf{e}_i)=\varphi(p^{\lambda_i}\mathbf{e}_i)=\varphi(F \mathbf{e}_i)=p^\ell\varphi(\mathbf{e}_i)^\sigma.\end{equation}
This is only possible if the $p$-adic valuation of the coefficients on both sides is the same, thus $\lambda_i=\ell$ for all $i$. This shows that at each point the $\sigma$-module $M$ has slope $\ell\in\NN$. 

By Dieudonn\'e--Manin, $M$ is locally at a point $\mathfrak{p}$ isogenous to $E(\ell)^n$ for some $n\in\NN$. Since $\ell$ is integral, the $E(\ell)$ are of rank one as $W(\overline{\kappa})$-modules, hence $n=\rank(M)$. Since $W^\dagger(\overline{\kappa})=W(\overline{\kappa})$, the images of the $\mathbf{e}_i$ under (the localisation of) $\varphi$ are in fact in $K$, thus cannot by linearly independent over $K$ --- and neither can the $\mathbf{e}_i$ by injectivity of $\varphi$, unless $n=1$.
\end{proof}

\begin{lem}
The preimage of the overconvergent Witt ring $\varphi^{-1}\left(W^\dagger(\overline{A})\right)$ is non-trivial.
\end{lem}

\begin{proof}
Here we use the fact, that the convergent and overconvergent Witt vectors over a ring $\overline{A}$ give rise to \'etale sheaves over $\Spec\overline{A}$. By the previous lemma we know that locally around a point $\mathfrak{p}$ the $\sigma$-module $M$ is of rank one, more precisely, up to isogeny of the form $E(\ell)$ and we can choose a local section $\mathbf{e}$ as above that generates it, so that it satisfies the equation
$$p^{\ell}(\varphi(e))=\varphi(p^{\ell}e)=\varphi(F(e)) =p^{\ell}\varphi(e)^\sigma$$
thus $\varphi(e)$ is contained in the kernel of $1-\sigma$, thus in $W(\overline{\kappa})$ at $\mathfrak{p}$. After multiplication with a power of $p$, in order to make up for the isogeny, we see that in an \'etale neighbourhood of $\mathfrak{p}$, the linearity of $\varphi$ with respect to $W^\dagger(\overline{A})$ shows that the preimage of overconvergent elements in this neighbourhood is nontrivial. We can do this for each point $\mathfrak{p}$, but in fact it is sufficient to look at a finite (fine enough) cover of $\Spec\overline{A}$, so we only have to multiply with a finite power of $p$. 
\end{proof}

Now we are ready to prove the main point of this section. 

\begin{theo}\label{LocalMain}
With $M_1$ and $M_2$ overconvergent $(\sigma,\nabla)$-modules as above, let $f:M_1\otimes W(\overline{A})\rightarrow M_2\otimes W(\overline{A})$ be a morphism of $(\sigma,\nabla)$-modules. Then there exists a morphism $g:M_1\rightarrow M_2$ of overconvergent $(\sigma,\nabla)$-modules that induces $f$.
\end{theo}

\begin{proof}
It suffices to show the claim for a morphism $\varphi:M\rightarrow W(\overline{A})$, where $M$ is an overconvergent $(\sigma,\nabla)$-module, and $\varphi$ a morphism of convergent $(\sigma,\nabla)$-modules that is $W^\dagger(\overline{A})$-linear, and satisfies for some $\ell\in\NN$ the conditions
\begin{enumerate}\item for all $\mathbf{e}\in M$, 
$$\varphi(F\mathbf{e})=p^\ell\varphi(\mathbf{e})$$
\item for all $\mathbf{e}\in M$,
$$\varphi(\nabla\mathbf{e})=d\varphi(\mathbf{e}).$$
\end{enumerate}
Although the dual $M^\ast_1:=\Hom(M_1,W^\dagger(\overline{A}))$ is not necessarily a $(\sigma,\nabla)$-module --- only after inverting $p$ (compare \cite[p.5]{Kedlaya4}) --- its Tate twist $M_1^\ast(\ell)$ for a sufficiently large $\ell$ is. Furthermore, there is a canonical isomorphism $\Hom(M_1,M_2)(\ell)\cong M^\ast_1(\ell)\otimes M_2$ of overconvergent $(\sigma,\nabla)$-modules. Via this isomorphism, $f$ induces a $W^\dagger(\overline{A})$-linear morphism $\varphi:M:=M_1^\ast(\ell)\otimes M_2\rightarrow W(\overline{A})$, compatible with the differential graded complexes, such that the above conditions are satisfied. 

Accordingly the task now is to show that such a $\varphi$ has image in $W^\dagger(\overline{A})$. 

Moreover, we can assume that $\varphi$ is injective. Otherwise, let $H\subset  M$ be the kernel of $\varphi$, which is obviously closed under $\sigma$ and $\nabla$, so that the quotient $M\slash H$ is again a $(\sigma,\nabla)$-module, and the induced map 
$$\overline{\varphi}: M\slash H\rightarrow W(\overline{A})$$
injective. With the same argument, we can therefore assume, that $M$ is of rank 1. 

Let $P$ be the preimage of $W^\dagger(\overline{A})$ under $\varphi$. As seen in the previous lemma, $P$ is non-trivial, and what is more, of full rank 1 and slope $\ell$. Hence it is clear (by Dieudonn\'e--Manin again) that $P=M$ and trivially closed under $\nabla$. We conclude therefore that $\varphi(M)=W^\dagger(\overline{A})$.
\end{proof}

\section{Global full faithfulness}

In this section we show how to reduce global full faithfulness of the forgetful functor 
$$j^\ast: F\text{-}\Con^\dagger(X/W(k))\rightarrow F\text{-}\Con(X/W(k))$$ 
to the local case, thereby proving the main result. The discussion borrows an argument from  \cite[Sec. 7]{Kedlaya4}.

In the basic set-up, that is, $k$ a perfect field of characteristic $p\neq 0$, $W(k)$ its ring of Witt vectors, $X$ a smooth variety over $k$, consider two overconvergent $(\sigma,\nabla)$-modules $M_1, M_2 \in F\text{-}\Con^\dagger(X/W(k))$. We have to show that the set of morphisms $M_1\rightarrow M_2$ in $F\text{-}\Con^\dagger(X/W(k))$ and $F\text{-}\Con(X/W(k))$ are equivalent. 

The set $\underline{\Hom}(M_1,M_2)=M$ is again an overconvergent $(\sigma,\nabla)$-module. The morphisms from $M_1$ to $M_2$ in the overconvergent category $F\text{-}\Con^\dagger(X/W(k))$ as well as in the convergent category $F\text{-}\Con(X/W(k))$ correspond to the elements of $M$ that are killed by $\nabla$ and fixed by $F$. Let's denote this set by $\h^0_{F}(X,M)$. Hence we can assert the following.

\begin{lem}\label{LemmaSuffices}
To show full faithfulness of $j^\ast: F\text{-}\Con^\dagger(X/W(k))\rightarrow F\text{-}\Con(X/W(k))$, it suffices to show that for $M\in F\text{-}\Con^\dagger(X/W(k))$ the rank of $\h^0_{F}(X,M)$ is the same in the convergent and overconvergent category.
\end{lem}

The definition of $\h^0_{F}(X,M)\subset \h^0(X,M)=\HH^0(M\otimes W^\dagger\Omega_{X/k})$ is local. Therefore, one only has to choose a suitable cover of $X$ in order to be able to apply the theorem of the previous section. But this is possible by the following result taken from Christopher Davis' thesis.

\begin{lem}[Davis 2009] 
Any smooth scheme has an open cover by standard \'etale affines. Moreover, we can choose this cover to consist of sets so that any finite intersection is again standard \'etale affine.
\end{lem}

Here standard \'etale affine means a ring of the form $\overline{B} = (\overline{A}_h[z]/P(z))_g$
with $\overline{A}$ a polynomial algebra  $h\in \overline{A}$, $\overline{A}_h$ the localisation, $P(z)$ monic, $g\in \overline{A}_h[z]/P(z)$, $(\overline{A}_h[z]/P(z))_g$ again the localisation, and $P'(z)$ invertible in $\overline{B}$.

Thus it is enough to prove what is suggested in Lemma \ref{LemmaSuffices} for $X$ such a standard \'etale affine, in which case this is equivalent to what was shown in Theorem \ref{LocalMain} in the previous section.

\begin{rem}
The same result holds true if we replace $F$ by some power $F^a$, with $a\in\NN$.
\end{rem}



\addcontentsline{toc}{section}{References}

\vspace{1cm}
\hrule
\vspace{.5cm}
{\noindent
\textsc{Universit\"at Regensburg}\\     	
Fakult\"at f\"ur Mathematik	\\
Universit\"atsstra{\ss}e 31	\\
93053 Regensburg		\\
Germany		      	\\
(+ 49) 941-943-2664       	\\
\verb=veronika.ertl@mathematik.uni-regensburg.de=\\
\url{http://www.mathematik.uni-regensburg.de/ertl/} \\}

\end{document}